\documentclass[11pt]{article}
\usepackage{amsfonts}
\usepackage{latexsym, amssymb, amsmath, amscd, amsfonts, epsfig, graphicx, colordvi}
\usepackage{ifpdf}
\usepackage{graphicx}
\usepackage{xcolor}
\usepackage[dvipsnames]{pstricks}
\newtheorem{thm}{Theorem}[section]

\newtheorem{defi}[thm]{Definition}
\newtheorem{lem}[thm]{Lemma}

\def\pf{\noindent{\it Proof.} }
\def\qed{\nopagebreak\hfill{\rule{4pt}{7pt}}
\medbreak}

\def\qed{\nopagebreak\hfill{\rule{4pt}{7pt}}
\medbreak}

\makeatletter
\def\ExtendSymbol#1#2#3#4#5{\ext@arrow 0099{\arrowfill@#1#2#3}{#4}{#5}}
\makeatother

\title{Parity Considerations  in Rogers--Ramanujan--Gordon Type Overpartitions }
\author{Doris D. M. Sang\raisebox{5pt}{\scriptsize 1} and Diane Y. H.
Shi\raisebox{5pt}{\scriptsize 2}}
\date{
\vspace{15pt}\raisebox{5pt}{\scriptsize
1\,}School of Mathematics and Quantitative
Economics\\Dongbei
University of Finance and Economics, Liaoning 116025, P.R. China\\
\vspace{15pt}\raisebox{5pt}{\scriptsize
2\,}School of
Mathematics, \\ Tianjin University, Tianjin 300072, P.R. China\vskip 0.2 cm
Email: $^1$sangdm@dufe.edu.cn, $^2$shiyahui@tju.edu.cn}

\begin{document}
\bibliographystyle{plain} %
\maketitle
\noindent {\bf Abstract.}
In 2010, Andrews considers a variety of parity questions connected to
classical partition identities of Euler, Rogers, Ramanujan and Gordon.  As a large part in his paper, Andrews considered the partitions by restricting the parity of occurrences of even numbers or  odd numbers  in the Rogers-Ramanujan-Gordon type.
The Rogers--Ramanujan--Gordon type partition was defined by Gordon in 1961 as a combinatorial generalization of the Rogers--Ramaujan identities with odd moduli.
 In 1974, Andrews derived an identity which can be considered as the generating function counterpart of
the Rogers--Ramanujan--Gordon theorem, and since then it has been called the Andrews--Gordon identity. By revisting the Andrews--Gordon identity Andrews extended
his results by considering some additional restrictions involving parities to obtain  some Rogers--Ramanujan--Gordon type theorems and Andrews--Gordon type
identities. In the end of Andrews' paper, he posed $15$ open problems. Most of Andrews' $15$ open problems have been settled, but the $11$th that ``extend the parity
indices to overpartitions in a manner" has not. In 2013, Chen, Sang and Shi, derived the overpartition analogues of the Rogers--Ramanujan--Gordon theorem and the  Andrews--Gordon identity.  In this paper, we post some parity restrictions on these overpartitions analogues to get some Rogers--Ramanujan--Gordon type overpartition theorems.



\section{ Introduction}

The celebrated combinatorial generalization of the Rogers--Ramanujan identities was given by Gordon \cite{gor61} in 1961, which can be  stated as follows:

\begin{thm}(Rogers-Ramanujan-Gordon)\label{Gordon} For $k\geq a\geq 1$,
let $B_{k,a}(n)$ denote the number of partitions of $n$ of the form
$\lambda_1 + \lambda_2 + \cdots + \lambda_s$, where $\lambda_j \geq
\lambda_{j+1}$, $\lambda_j-\lambda_{j+k-1}\geq2$ and part $1$
appears at most $a-1$ times.
 Let $A_{k,a}(n)$ denote the number of
partitions of $n$ into parts $\not \equiv0,\pm a\; \pmod{2k + 1}$.
Then for any $n \geq0$, we have
 \begin{equation*}
A_{k,a}(n) = B_{k,a}(n).
 \end{equation*}
\end{thm}

In 1967, Andrews \cite{and74} established his analytic generalization of the Rogers--Ramanujan identities with the moduli from $5$ to all odd positive integers:
\begin{thm}For $k\geq a\geq 1$, we have
 \begin{equation}
 \label{eqRRG}\sum_{N_1\geq N_2\geq\cdots\geq
N_{k-1}\geq0}\frac{q^{N_1^2+N_2^2+\cdots+N_{k-1}^2+N_{a+1}+\cdots+N_{k-1}}
}{(q)_{N_1-N_2}\cdots(q)_{N_{k-2}-N_{k-1}}(q)_{N_{k-1}}}=\frac{(q^a,q^{2k+1-a},
q^{2k+1};q^{2k+1})_\infty}{(q)_\infty}.
\end{equation}
\end{thm}
Here and in the rest of this paper, we  adopt  the common notation  as used in Andrews \cite{and76}. Let
 \[(a)_\infty=(a;q)_\infty=\prod_{i=0}^{\infty}(1-aq^i),\]
 and \[(a)_n=(a;q)_n=\frac{(a)_\infty}{(aq^n)_\infty}.\]
We also write
 \[(a_1,\ldots,a_k;q)_\infty=(a_1;q)_\infty\cdots(a_k;q)_\infty.\]

Let us give an overview of  some definitions. A partition $\lambda$
of a positive integer $n$ is a non-increasing sequence of positive
integers $\lambda_1\geq \cdots\geq \lambda_s>0$ such that
$n=\lambda_1+\cdots+\lambda_s$. The partition of zero is the
partition with no parts. An overpartition $\lambda$ of a positive
integer $n$ is also a non-increasing sequence of positive integers
$\lambda_1\geq \cdots\geq \lambda_s>0$ such that
$n=\lambda_1+\cdots+\lambda_s$ and the first occurrence  of  each
integer may be overlined. For example,
$(\overline{7},7,6,\overline{5},2,\overline{1})$ is an overpartition
of $28$. For a partition
or an overpartition $\lambda$ and for any integer $l$, let
$f_l(\lambda) (f_{\overline{l}}(\lambda))$ denote the number of
occurrences of $l$ non-overlined (overlined) in $\lambda$.

 In 2010, Andrews \cite{and10} investigated a variety of parity questions in classical partition identities. In particular, he  revisited his generating function \eqref{eqRRG} and extended
his results by considering some additional restrictions involving parities. The first theorem considered that in the case  $k\equiv a (\rm{mod}\ 2)$, the conditions are subject to that the occurrences of even numbers must be even.
\begin{thm} Suppose $k\geq a\geq 1$ are integers with $k\equiv a \pmod{2}$. Let $W_{k,a}(n)$
denote the number of those partitions enumerated by $B_{k,a}(n)$ with the added restriction
that even parts appear an even number of times. If $k$ and $a$ are both even, let
$G_{k,a}(n)$ denote the number of partitions of $n$ in which no odd part is repeated and
no even part $\equiv0\pm a \pmod{2k + 2}$. If $k$ and $a$ are both odd, let $G_  {k,a}(n)$ denote
the number of partitions of $n$ into parts that are neither $\equiv 2 \pmod{4}$ nor $\equiv0\pm a
\pmod{2k + 2}$. Then for all $n\geq 0$
\[ W_{k,a}(n) = G_{k,a}(n).\]
\end{thm}
The generating function form of this theorem can be state as follows:
\begin{thm}For $k\geq a\geq 1$ and $k\equiv a \pmod{2}$, we have
\begin{align}\nonumber\sum_{n\geq 0}W_{k,a}(n)q^n&=\sum_{N_1\geq N_2\geq\cdots\geq
N_{k-1}\geq0}\frac{q^{N_1^2+N_2^2+\cdots+N_{k-1}^2+2N_{a}+2N_{a+2}+\cdots+2N_{k-2}}
}{(q^2;q^2)_{N_1-N_2}\cdots(q^2;q^2)_{N_{k-2}-N_{k-1}}(q^2;q^2)_{N_{k-1}}}
\\[6pt]\label{WG}&\qquad=\frac{(-q;q^2)_\infty(q^a,q^{2k+2-a},q^{2k+2};q^{2k+2})_\infty}{(q^2;q^2)_\infty}.
\end{align}
 \end{thm}

When  $k\not\equiv a \pmod{2}$, Andrews  considered the partitions provided that the occurrences of odd parts must even.
\begin{thm}
 Suppose $k\geq a\geq 1$ with $k$ odd and $a$ even. Let $\overline{W}_{k,a}(n)$ denote the
number of those partitions enumerated by $B_{k,a}(n)$ with added restriction that odd
parts appear an even number of times. Then
 \[\sum_{n\geq0}\overline{W}_{k,a}(n)q^n=\frac{(q^a,q^{2k+2-a},q^{2k+2};
 q^{2k+2})_\infty}{(-q;q^2)_\infty(q;q)_\infty}.\]
\end{thm}
The generating function form  of this theorem is also given by Andrews \cite{and10} as follows:
\begin{thm}For  $k\geq a\geq 1$ with $k$  odd and $a$ even, we have
\begin{align}\nonumber&\sum_{N_1\geq N_2\geq\cdots\geq
N_{k-1}\geq0}\frac{q^{N_1^2+N_2^2+\cdots+N_{k-1}^2+n_1+n_3+\cdots+n_{a-3}+N_{a-1}
+N_a+\cdots+N_{k-1}}
}{(q^2;q^2)_{N_1-N_2}\cdots(q^2;q^2)_{N_{k-2}-N_{k-1}}(q^2;q^2)_{N_{k-1}}}
\\\nonumber&\quad=\frac{(-q^2;q^2)_\infty(q^a,q^{2k+2-a},q^{2k+2};q^{2k+2})_\infty}
{(q^2;q^2)_\infty}\\&\quad=\frac{(q^a,q^{2k+2-a},q^{2k+2};q^{2k+2})_\infty}
{(-q;q^2)_\infty(q;q)_\infty}.
\end{align}
\end{thm}
Those generalizations are obtained by using double recurrences  and the Defining
$q$-Difference Equations Principle.

Andrews just considered  $W_{k,a}(n)$ in the case that $k\equiv a (\rm{mod}\ 2)$ and  $\overline{W}_{k,a}(n)$ in the case that $k$  is odd and $a$  is even. In 2013, Kim and Yee \cite{kim13} derived the generating function of $W_{k,a}(n)$ and $\overline{W}_{k,a}(n)$
with $k$ and $a$ in other parities.
For $k$ and $a$ have different parities,  Kim and Yee derived the following result about $W_{k,a}(n)$.
\begin{thm}For $k\geq a\geq1$, $k\not\equiv a\pmod{2}$,
\begin{align}\nonumber\sum_{n\geq 0}W_{k,a}(n)q^n=&\frac{(-q^3;q^2)_\infty(q^{a+1},q^{2k+1-a},q^{2k+2};q^{2k+2})_\infty}{(q^2;q^2)_\infty}
\\&+\frac{q(-q^3;q^2)_\infty(q^{a-1},q^{2k+3-a},q^{2k+2};q^{2k+2})_\infty}{(q^2;q^2)_\infty}.
\end{align}
\end{thm}
For $\overline{W}_{k,a}(n)$, Kim and Yee have the following result which is a ``missing" case of Andrews.
\begin{thm}
Suppose $k\geq a\geq 1$ with $k$ even and $a$ odd. Then
\begin{align}
\sum_{n\geq 0}\overline{W}_{k,a}(n)q^n=\frac{(-q^2;q^2)_\infty(q^{a+1},q^{2k+1-a},q^{2k+1};q^{2k+1})_\infty}{(q^2;q^2)_\infty}.
\end{align}
\end{thm}
They also note the following relation for $\overline{W}_{k,a}(n)$.
\begin{thm}For $k\geq a\geq1$ with $a$ even and $n\geq 1$.
\begin{equation}\overline{W}_{k,a}(n)=\overline{W}_{k,a-1}(n).
\end{equation}
\end{thm}
Then the generating function of $W_{k,a}(n)$ and $\overline{W}_{k,a}(n)$ in all parities have been derived.
 In this paper, we will give parity restrictions on the Rogers--Ramanujan--Gordon type overpartitions. We first revisit the Rogers--Ramanujan--Gordon type overpartition theorems.

In 2012, Chen, Sang and Shi \cite{chen} established the overpartition analogue of the Rogers--Ramanujan--Gordon theorem:
\begin{thm}[The Rogers--Ramanujan--Gordon type overpartition theorem]\label{thmlast1} For $k\geq 2$ and  $k\geq i\geq 1$, let $\overline{B}_{k,i}(n)$ denote
the number of overpartitions of $n$ of the form $\lambda_1+\lambda_2+\cdots+\lambda_s$
such that part $1$  occurs as a non-overlined part at most $i-1$ times,
and  $\lambda_j-\lambda_{j+k-1}\geq1$ if $\lambda_j$ is overlined and
$\lambda_j-\lambda_{j+k-1}\geq2$ otherwise. For $k> i\geq 1$, let $\overline{A}_{k,i}(n)$ denote the number of
overpartitions of $n$ whose non-overlined parts are not congruent to
$0,\pm i$ modulo $2k$ and let $\overline{A}_{k,k}(n)$ denote the number of overpartitions of $n$ with parts not divisible by $k$.
 Then we have
 \begin{equation}
 \overline{A}_{k,i}(n)=\overline{B}_{k,i}(n).
 \end{equation}
\end{thm}

We shall focus on  the parity of the frequencies of a number $l$ by consider $f_{l}+f_{\overline{l}}$. We give the following two definitions.
\begin{defi}For  $k\geq a\geq 1$, let $U_{k,a}(n)$ denote
the number of overpartitions of $n$ of the form $(\overline{1}f_{\overline{1}},\ 1f_1,\ \overline{2}f_{\overline{2}},\ 2f_2, \ldots)$
such that
\begin{itemize}
\item[(i)]$f_1(\lambda)\leq a-1+f_{\overline{1}}(\lambda)$;
\item[(ii)]$f_{2l-1}(\lambda)\geq f_{\overline{2l-1}}(\lambda)$;
\item[(iii)]$f_{2l}(\lambda)+f_{\overline{2l}}(\lambda)\equiv 0\pmod{2}$;
\item[(iv)]$f_l(\lambda)+f_{\overline{l}}(\lambda)+f_{l+1}(\lambda)\leq k-1+f_{\overline{l+1}}(\lambda)$.
\end{itemize}
\end{defi}

\begin{defi}For  $k\geq a\geq 1$, let $\overline{U}_{k,a}(n)$ denote
the number of overpartitions of $n$ of the form $(\overline{1}f_{\overline{1}},\ 1f_1,\ \overline{2}f_{\overline{2}},\ 2f_2, \ldots)$
such that
\begin{itemize}
\item[1.]$f_1(\lambda)\leq a-1+f_{\overline{1}}(\lambda)$;
\item[2.]$f_{2l}(\lambda)\geq f_{\overline{2l}}(\lambda)$;
\item[3.]$f_{2l-1}(\lambda)+f_{\overline{2l-1}}(\lambda)\equiv 0 \pmod{2}$;
\item[4.]$f_l(\lambda)+f_{\overline{l}}(\lambda)+f_{l+1}(\lambda)\leq k-1+f_{\overline{l+1}}(\lambda)$.
\end{itemize}
\end{defi}
 Here are our main results.
\begin{thm}\label{u22}For  $k\geq a\geq 1$ and $k\equiv a \pmod{2}$,  we have
\begin{equation}\sum_{n\geq 0}U_{k,a}(n)q^n=\frac{(-q;q)_\infty(q^a,q^{2k-a},q^{2k};q^{2k})_\infty}{(q^2;q^2)_\infty}.
\end{equation}
\end{thm}
\begin{thm}\label{u21}For  $k\geq a\geq 1$ and $k\not\equiv a\pmod{2}$, we have
\begin{equation}\sum_{n\geq0}U_{k,a}(n)q^n=\frac{(-q^2;q)_\infty(q^{a+1},q^{2k-a-1},q^{2k};q^{2k})_\infty}{(q^2;q^2)_\infty}
+\frac{xq(-q^2;q)_\infty(q^{a-1},q^{2k-a+1},q^{2k};q^{2k})_\infty}{(q^2;q^2)_\infty}.
\end{equation}
\end{thm}

\begin{thm}\label{overlineu22}For  $k\geq a\geq 2$ with $a$ even, we have
\begin{equation}\sum_{n\geq 0}\overline{U}_{k,a}(n)q^n=\sum_{n\geq 0}\overline{U}_{k,a-1}(n)q^n=\frac{(-q^2;q^2)^2_\infty(q^{a},q^{2k-a},q^{2k};q^{2k})_\infty}{(q^{2};q^2)_\infty}.
 \end{equation}
\end{thm}

By Theorem \ref{u22}--Theorem \ref{overlineu22}, we get the generating functions of $U_{k,a}(n)$ and $\overline{U}_{k,a}(n)$ in all parities of $k$ and $a$.

This paper is organized as follows. In Section 2, we shall recall $\overline{Q}_{k,a}(x;q)$ which has been used in the proof of Theorem \ref{thmlast1}. And recall the defining
$q$-difference equations principle.
In Section 3, we prove Theorem \ref{u22} in the case that $k\equiv a\equiv 0\pmod{2}$. In Section 4, we derive the generating function of $U_{k,a}$ in the case that $k\not \equiv a \pmod{2}$. In Section 5, we give the generating function of $\overline{U}_{2k,2a}(n)$. In Section 6, we derive the generating function of $U_{2k+1,2a+1}(n)$ and $\overline{U}_{2k+1,2a+1}(n)$ by two $q$-differential relations. In Section 7, we give the combinatorial interpretations of two identities, which reveal the relations among $U_{k,a}(n)$, $\overline{U}_{k,a}(n)$ and $\overline{B}_{k,a}(n)$.

\section{Background}
We shall use a series $H_{k,i}(a;x;q)$  introduced by
 Andrews \cite{and66, and74}, which is defined by
\begin{equation}\label{eqH}
H_{k,i}(a;x;q)=
\sum_{n=0}^{\infty}\frac{x^{kn}q^{kn^2+n-in}a^n(1-x^iq^{2ni})
(axq^{n+1})_{\infty}(1/a)_n}{(q)_n(xq^n)_\infty}.
\end{equation}
In his algebraic proof of the Rogers--Ramanujan--Gordon theorem,
Andrews used the function $J_{k,i}(a;x;q)$ constructed based on $H_{k,i}(a;x;q)$,
\begin{equation*} 
J_{k,i}(a;x;q)=
H_{k,i}(a;xq;q)-axqH_{k,i-1}(a;xq;q).
\end{equation*}
Chen, Sang and Shi \cite{chen} proved Theorem \ref{thmlast1} by using a specialization of $H_{k,i}(a;x;q)$ which we denote it here by $\overline{Q}_{k,i}(x;q)$, that 
\begin{equation}\label{eqH}
\overline{Q}_{k,i}(x;q) =H_{k,i}(-1/q;xq;q)
=\sum_{n=0}^{\infty}\frac{(-1)^nx^{kn}q^{kn^2+kn-in}(1-x^iq^{(2n+1)i})
(-xq^{n+1})_{\infty}(-q)_n}{(q)_n(xq^{n+1})_\infty},
\end{equation}
 to get the following result
\[\overline{Q}_{k,i}(x;q)=\sum_{m,n\geq 0}\overline{B}_{k,i}(n)x^mq^n.\]

In this paper, we also employ $\overline{Q}_{k,i}(x;q)$ to prove Theorem \ref{u22}--Theorem \ref{overlineu22}. We stated some results on $\overline{Q}_{k,i}(x;q)$ which is derived from Andrews' relations about $H_{k,i}(a;x;q)$ and have been stated by Chen, Sang and Shi in \cite{chen}.
\begin{lem}\label{lemQini}Let $\overline{Q}_{k,i}(x;q)=\sum_{m,n\geq 0}\overline{Q}_{k,i}x^mq^n$, then we have the following initial values
\begin{equation}\overline{Q}_{k,i}(0,0)=1, \text{ for\ } k\geq i\geq 1,\end{equation}
\begin{equation}\label{Qk0}\overline{Q}_{k,0}(m,n)=0, \text{ for } k\geq 1,\ m,n\geq 0,\end{equation}
\begin{equation}\overline{Q}_{k,i}(m,n)=0, \text{ if }m\text{ or }n\text{ is zero but not both},\end{equation}
\begin{equation}\overline{Q}_{k,i}(m,n)=0, \text{ if }m\text{ or }n<0.\end{equation}\end{lem}
By \eqref{eqH}, one can derive the following relation.
\begin{lem}
For $k\geq 1$, we have 
\begin{equation}\label{Qminus}\overline{Q}_{k,-\frac{1}{2}}(x;q)=-(xq)^{-\frac{1}{2}}\overline{Q}_{k,\frac{1}{2}}(x;q)
\end{equation}\end{lem}
\begin{lem}\cite[P106,Lemma 7.1]{and76}\label{lemQ}We have
\begin{equation}\label{eqJ}
\overline{Q}_{k,i}(x;q)-\overline{Q}_{k,i-1}(x;q)=(xq)^i\overline{Q}_{k,k-i}(xq;q)+(xq)^{i-1}\overline{Q}_{k,k-i+1}(xq;q).
\end{equation}
\end{lem}
By relation \eqref{eqJ}, we shall give the following relation satisfied by $\overline{Q}_{k,i}(m,n)$.
\begin{thm}\label{Qka}For $k\geq a\geq 1$, we have
\begin{align}
  \nonumber\overline{Q}_{k,a}(x;q)&=\sum_{i=1}^a (xq)^i \left[\sum_{h=1}^{k-i}(xq^2)^h\overline{Q}_{k,k-h}(xq^2;q) +\sum_{h=0}^{k-i-1}(xq^2)^h\overline{Q}_{k,k-h}(xq^2;q)\right]
  \\&\quad+ \label{QKI} \sum_{i=0}^{a-1}(xq)^i\left[\sum_{h=1}^{k-i}(xq^2)^h\overline{Q}_{k,k-h}(xq^2;q)+\sum_{h=0}^{k-i-1}(xq^2)^h\overline{Q}_{k,k-h}(xq^2;q)\right].
  \end{align}
\end{thm}
\pf By relation \eqref{eqJ}, we have that
\begin{equation*}
  \overline{Q}_{k,a}(x;q)=\overline{Q}_{k,a-1}(x;q)+(xq)^a\overline{Q}_{k,k-a}(xq;q)+(xq)^{a-1}\overline{Q}_{k,k-a+1}(xq;q).\end{equation*}
  Then by successively using \eqref{Qka} to $\overline{Q}_{k,a-1}(x;q)$ in above identity, we have that \begin{align*}\overline{Q}_{k,a}(x;q)&=\overline{Q}_{k,a-2}(x;q)+(xq)^{a-1}\overline{Q}_{k,k-a+1}(xq;q)+(xq)^{a-2}\overline{Q}_{k,k-a+2}(xq;q)
  \\&\qquad+(xq)^a\overline{Q}_{k,k-a}(xq;q)+(xq)^{a-1}\overline{Q}_{k,k-a+1}(xq;q)
  \\&\quad=\overline{Q}_{k,0}(x;q)+\sum_{i=1}^a(xq)^i\overline{Q}_{k,k-i}(xq;q)+\sum_{i=0}^{a-1}(xq)^i\overline{Q}_{k,k-i}(xq;q).
  \end{align*}
  By \eqref{Qk0}, we have
  \begin{equation*}\label{Qka2}\overline{Q}_{k,a}(x;q)=\sum_{i=1}^a(xq)^i\overline{Q}_{k,k-i}(xq;q)+\sum_{i=0}^{a-1}(xq)^i\overline{Q}_{k,k-i}(xq;q).\end{equation*}
  Now we successively employ  \eqref{Qka2} to $\overline{Q}_{k,k-i}(xq;q)$, to get that
  \begin{align*}
  &\overline{Q}_{k,a}(x;q)=\sum_{i=1}^a (xq)^i \left[\sum_{h=1}^{k-i}(xq^2)^h\overline{Q}_{k,k-h}(xq^2;q) +\sum_{h=0}^{k-i-1}(xq^2)^h\overline{Q}_{k,k-h}(xq^2;q)\right]
  \\&\quad+  \sum_{i=0}^{a-1}(xq)^i\left[\sum_{h=1}^{k-i}(xq^2)\overline{Q}_{k,k-h}(xq^2;q)+\sum_{h=0}^{k-i-1}\overline{Q}_{k,k-h}(xq^2;q)\right].
  \end{align*}\qed
Recall the Defining
$q$-Difference Equations Principle as follows \cite{and10}.
\begin{thm}Suppose that we  have two sets of functions, for $1\leq i\leq r$, $f_i(x,q)$ and $g_i(x,q)$,  which are analytic in $x$ and $q$ for $|q|<1$ and $|x|<|q|^{-1}$. Further,
suppose that for each $i$, $f_i(0,q)=g_i(0,q)$,
\begin{equation}
f_i(x,q)=\sum_{j=1}^rh_{i,j}(x,q)f_j(xq^{e(i,j)},q),
\end{equation}
and
\begin{equation}
g_i(x,q)=\sum_{j=1}^rh_{i,j}(x,q)g_j(xq^{e(i,j)},q),
\end{equation}
where the $e(i,j)$ are all positive integers and the $h_{i,j}(x;q)$ are polynomials in $x$ and
$q$. Then it follows by a double induction on the double power series coefficients that
\[f_i(x,q)=g_i(x,q), \qquad 1\leq i\leq r.\]
\end{thm}

We can see that by the defining
$q$-difference equations principle, Lemma \ref{lemQini} and  Theorem \ref{Qka}, the functions $\overline{Q}_{k,a}(x;q)$ are uniquely determined.

\section{The generating function of $U_{2k,2a}(n)$}

In this section and in the sequel, we let $U_{k,a}(m,n)$ (resp. $\overline{U}_{k,a}(m,n)$) denote the number of overpartitions enumerated by $U_{k,a}(n)$ (resp. $\overline{U}_{k,a}(n)$) that have exactly $m$ parts. The
related generating functions are
\[U_{k,a}(x;q)=\sum_{m,n\geq 0}U_{k,a}(m,n)x^mq^n,\]
and
\[\overline{U}_{k,a}(x;q)=\sum_{m,n\geq 0}\overline{U}_{k,a}(m,n)x^mq^n.\]
In this section, we shall prove the even case of  Theorem \ref{u22} i.e. $k\equiv a\equiv0\pmod{2}$ from  following result
\begin{thm}\label{UQ} For $k\geq a\geq 1$ are integers, we have
\begin{equation}\label{eqUQ}U_{2k,2a}(x;q)=(-xq;q^2)_\infty\overline{Q}_{k,a}(x^2;q^2).
\end{equation}
\end{thm}

In previous  section, we have stated some relations satisfied by $\overline{Q}_{k,a}(x;q)$. Now we define
\[\overline{T}_{2k,2a}(x;q):=(-xq;q^2)_\infty\overline{Q}_{k,a}(x^2;q^2),\]
then shall prove that $\overline{T}_{2k,2a}(x;q)$ and $U_{2k,2a}(x;q)$ satisfy the same $q$-differential relations with the same initial values, implying \[U_{2k,2a}(x;q)=\overline{T}_{2k,2a}(x;q).\]
By the definition of $\overline{T}_{2k,2a}(x;q)$ and \eqref{QKI}, we have.
\begin{align}\label{Tka}
  \nonumber\overline{T}_{2k,2a}(x;q)&=(1+xq)\sum_{i=1}^a (xq)^{2i} \left[\sum_{h=1}^{k-i}(xq^2)^h\overline{T}_{2k,2k-2h}(xq^2;q) +\sum_{h=0}^{k-i-1}(xq^2)^h\overline{T}_{2k,2k-2h}(xq^2;q)\right]
  \\&\quad+ (1+xq) \sum_{i=0}^{a-1}(xq)^{2i}\left[\sum_{h=1}^{k-i}(xq^2)^h\overline{T}_{2k,2k-2h}(xq^2;q)+\sum_{h=0}^{k-i-1}(xq^2)^h\overline{T}_{2k,2k-2h}(xq^2;q)\right].
  \end{align}

We shall derive  the $q$-differential relation satisfied by $U_{2k,2a}(x;q)$ by using the combinatorial tool.

\begin{thm}For $k\geq a\geq1$, we have \begin{align}\nonumber&U_{2k,2a}(x;q)-U_{2k,2a-2}(x;q)
\\\nonumber&=(1+xq)\left[(xq)^{2a}\sum_{h=1}^{k-a}(xq^2)^{2h}U_{2k,2k-2h}(xq^2;q)+(xq)^{2a}\sum_{h=0}^{k-a-1}(xq^2)^{2h}U_{2k,2k-2h}(xq^2;q)\right]
\\&+(1+xq)\left[(xq)^{2a-2}\sum_{h=1}^{k-a+1}(xq^2)^{2h}U_{2k,2k-2h}(xq^2;q)+(xq)^{2a-2}\sum_{h=0}^{k-a}(xq^2)^{2h}U_{2k,2k-2h}(xq^2;q)\right].\end{align}
\end{thm}
\pf By the definition of the overpartitions enumerated by $U_{2k,2a}(n)$ and $U_{2k,2a-2}(n)$ we know that $U_{2k,2a}(n)-U_{2k,2a-2}(n)$ is the number of overpartitions that enumerated by $U_{2k,2a}(n)$ with the restrictions that
\begin{itemize}
\item[1.] if there is an $\overline{1}$, then the number of $1$ is $2a$ or $2a-1$;
\item[2.] if there is no $\overline{1}$, then the number of $1$ is $2a-1$ or $2a-2$.
\end{itemize}

Then we consider the set in the following  four cases.

(i) The number of overpartitions enumerated by $U_{2k,2a}(n)$ with the restriction that is an $\overline{1}$ and the number of $1$ is $2a$.
If there is an $\overline{2}$, then the number of non-overlined $2$ is at most $2k-2a-1$ with the restriction that $f_{2l}+f_{\overline{2l}}$ is even.
Suppose that the number of $f_{2l}+f_{\overline{2l}}$ is $2h$, then the number of non-overlined $3$ is at most $2k-2h$ with one $\overline{3}$, and $2k-2h-1$ with no $\overline{3}$. We delete all the $1$'s and $2$'s, and then subtract $2$ from each other parts.  We get the overpartitions enumerated by $U_{2k,2k-2h}(m-2a-1-2h,n-2m+a+1)$. Then the generating function is
\begin{equation}\label{1}(xq)^{2a+1}\sum_{h=1}^{k-a}(xq^2)^{2h}U_{2k,2k-2h}(xq^2;q).\end{equation}

If there is no $\overline{2}$, then the number of non-overlined $2$ is at most $2k-2a-2$ with the restriction that $f_{2l}+f_{\overline{2l}}$ is even. Suppose that the number of $f_{2l}$ is $2h$, then the number of nonoverlined $3$ is at most $2k-2h$ with one $\overline{3}$ and $2k-2h-1$ with no $\overline{3}$. We delete all the $1$'s and $2$'s, and then subtract $2$ from each other parts.  We get the overpartitions enumerated by $U_{2k,2k-2h}(m-2a-1-2h,n-2m+2a+1)$. Then the generating function is \begin{equation}\label{2}(xq)^{2a+1}\sum_{h=0}^{k-a-1}(xq^2)^{2h}U_{2k,2k-2h}(xq^2;q).\end{equation}

(ii) Now, we consider the overpartitions that enumerated by $U_{2k,2a}(n)$ with the restriction that there is an $\overline{1}$ and the number of $1$ is $2a-1$. By similar analysis, we can get the generating function as
\begin{equation}\label{3}(xq)^{2a}\sum_{h=1}^{k-a}(xq^2)^{2h}U_{2k,2k-2h}(xq^2;q)+(xq)^{2a}\sum_{h=0}^{k-a-1}(xq^2)^{2h}U_{2k,2k-2h}(xq^2;q).\end{equation}

(iii) The number of overpartitions that enumerated by $U_{2k,2a}(n)$ with the restriction that there is no  $\overline{1}$ and the number of $1$ is $2a-1$.
The generating function is
\begin{equation}\label{4}(xq)^{2a-1}\sum_{h=1}^{k-a+1}(xq^2)^{2h}U_{2k,2k-2h}(xq^2;q)+(xq)^{2a-1}\sum_{h=0}^{k-a}(xq^2)^{2h}U_{2k,2k-2h}(xq^2;q).\end{equation}

(iv) The number of  overpartitions that enumerated by $U_{2k,2a}(n)$ with the restriction that there is an $\overline{1}$ and the number of $1$ is $2a-2$.
The generating function is
\begin{equation}\label{5}(xq)^{2a-2}\sum_{h=1}^{k-a+1}(xq^2)^{2h}U_{2k,2k-2h}(xq^2;q)+(xq)^{2a-2}\sum_{h=0}^{k-a}(xq^2)^{2h}U_{2k,2k-2h}(xq^2;q).\end{equation}

Compute the summation of \eqref{1}--\eqref{5} we complete the proof.
\qed

By successively employing the above theorem we get the following result.
\begin{thm}For $k\geq a\geq 1$, we have
\label{U22g}\begin{align}\nonumber&U_{2k,2a}(x;q)
\\\nonumber&=(1+xq)\bigg[\sum_{i=1}^a(xq)^{2i}\left(\sum_{h=1}^{k-i}(xq^2)^{2h}U_{2k,2k-2h}(xq^2;q)+\sum_{h=0}^{k-i-1}(xq^2)^{2h}U_{2k,2k-2h}(xq^2;q)\right)
\\&+\sum_{i=1}^{a}(xq)^{2i-2}\left(\sum_{h=1}^{k-i+1}(xq^2)^{2h}U_{2k,2k-2h}(xq^2;q)+\sum_{h=0}^{k-i}(xq^2)^{2h}U_{2k,2k-2h}(xq^2;q)\right)\bigg].
\end{align}
\end{thm}

\noindent \emph{The proof of Theorem \ref{UQ}}.
Now we  check the initial values.  One can easily get that  $U_{k,0}(m,n)=0$, $U_{k,a}(0,n)=U_{k,a}(m,0)=0$ and $U_{k,a}(0,0)=1$. Then $\overline{T}_{2k,2a}(x;q)$ and $U_{2k,2a}(x;q)$ have the same initial values. By Theorem \ref{U22g} and  the $q$-Difference Equations Principle \eqref{Tka}, we complete the proof.
\qed

\begin{thm}For $k\geq a\geq 1$, we have \begin{equation}\sum_{n\geq0}U_{2k,2a}(n)q^n=\frac{(-q;q)_\infty(q^{2a},q^{4k-2a},q^{4k};q^{4k})_\infty}{(q^2;q^2)_\infty}.
\end{equation}
\end{thm}
\pf Let $x=1$ in \eqref{eqUQ}, we have  \[U_{2k,2a}(1;q)=\sum_{m,n\geq0}U_{2k,2a}(m,n)q^n=\sum_{n\geq 0}U_{2k,2a}(n)q^n=\frac{(-q;q)_\infty(q^{2a},q^{4k-2a},q^{4k};q^{4k})_\infty}{(q^2;q^2)_\infty}.\]
\qed

\section{The generating function of $U_{2k,2a-1}(n)$ and $U_{2k+1,2a}(n)$}

In this section, we will give the generating function of $U_{2k,2a-1}(n)$ and $U_{2k+1,2a}(n)$.
For $1\leq i\leq a$, let $U^i_{k,a}(m,n)$ denote the number of overpartitions enumerated by $U_{k,a}(m,n)$ where if there is an $\overline{1}$, then the number of non-overlined $1$ is exactly $i$ otherwise  the number of non-overlined $1$ is exactly $i-1$, and the generating function denoted as follows.
\[U^{i}_{k,a}(x;q):=\sum_{m,n\geq 0}U^i_{k,a}x^mq^n.\]
\begin{thm}\label{21}For $k\geq a\geq 1$, we have
\begin{equation}\label{UU12}U_{2k,2a-1}(x;q)=U_{2k,2a}(x;q)-U^{2a}_{2k,2a}(x;q),
\end{equation}
and \begin{equation}\label{U22a}U^{2a}_{2k,2a}(x;q)=xq(-xq^3;q^2)_\infty[\overline{Q}_{k,a}(x^2;q^2)-\overline{Q}_{k,a-1}(x^2;q^2)].
\end{equation}
\end{thm}
\pf The relation \eqref{UU12} is a directly result by considering the number of $1$ and $\overline{1}$.

To derive \eqref{U22a}, we decompose the overpartions enumerated by $U^{2a}_{2k,2a}(m,n)$ into the following four cases.
\begin{itemize}
\item[1.] if there is an $\overline{1}$, then the number of nonoverlined  $1$ is $2a$. Suppose there is an $\overline{2}$  then the number of nonoverlined $2$ is odd and at most  $2k-2a-1$. Further, suppose there are exactly $2h-1$ non-overlined $2$, then by deleting all the $\overline{1}$, $1$'s, $\overline{2}$ and $2$'s, and subtracting $2$ from each of other parts, we get all overpartitions enumerated by $U_{2k,2k-2h}(m-2a-2h-1,n-2a-4h-1)$.
\item[2.] if there is an $\overline{1}$, then the number of nonoverlined  $1$ is $2a$. Futher, suppose there is no $\overline{2}$,  then the number of nonoverlined $2$ is odd and at most  $2k-2a-2$. Suppose there are exactly $2h$ non-overlined $2$, then by deleting all $\overline{1}$, $1$'s and $2$'s, and subtracting $2$ from each of other parts, we get all overpartitions enumerated by $U_{2k,2k-2h}(m-2a-2h-1,n-2a-4h-1)$.
 \item[3.] if there is no $\overline{1}$, then the number of nonoverlined  $1$ is $2a-1$. Suppose there is an $\overline{2}$,  then the number of nonoverlined $2$ is odd and at most  $2k-2a+1$. Further, suppose there are exactly $2h-1$ non-overlined $2$, then by deleting all $\overline{1}$, $1$'s, $\overline{2}$ and $2$'s, and subtracting $2$ from each of other parts, we get all overpartitions enumerated by $U_{2k,2k-2h}(m-2a-2h+1,n-2a-4h+1)$.
\item[4.] if there is no $\overline{1}$, then the number of nonoverlined  $1$ is $2a-1$. Suppose there is no $\overline{2}$  then the number of nonoverlined $2$ is odd and at most  $2k-2a$. Futher, suppose there are exactly $2h$ non-overlined $2$, then by deleting all $\overline{1}$, $1$'s, $\overline{2}$ and $2$'s, and subtracting $2$ from each of other parts, we get all overpartitions enumerated by $U_{2k,2k-2h}(m-2a-2h+1,n-2a-4h+1)$.
\end{itemize}
Thus, we get the following relation
 \begin{align*}
&U^{2a}_{2k,2a}(x;q)\\&=(xq)^{2a+1}\sum_{h=1}^{k-a}(xq^2)^{2h}U_{2k,2k-2h}(xq^2;q)
+(xq)^{2a+1}\sum_{h=0}^{k-a-1}(xq^2)^{2h}U_{2k,2k-2h}(xq^2;q)
\\&+(xq)^{2a-1}\sum_{h=1}^{k-a+1}(xq^2)^{2h}U_{2k,2k-2h}(xq^2;q)
+(xq)^{2a-1}\sum_{h=0}^{k-a}(xq^2)^{2h}U_{2k,2k-2h}(xq^2;q).
\end{align*}
In Section 3, we have proved that $U_{2k,2a}(x;q)=(-xq;q^2)_\infty\overline{Q}_{k,a}(x^2;q^2)$, then we have
\begin{align*}
&U^{2a}_{2k,2a}(x;q)\\&=(xq)^{2a+1}(-xq^3;q^2)_\infty\sum_{h=1}^{k-a}[(xq^2)^{2h}\overline{Q}_{k,k-h}(x^2q^4;q)+(xq^2)^{2h-2}\overline{Q}_{k,k-h+1}(x^2q^4;q)]
\\&+(xq)^{2a-1}(-xq^3;q^2)_\infty\sum_{h=1}^{k-a+1}[(xq^2)^{2h}\overline{Q}_{k,k-h}(x^2q^4;q)+(xq^2)^{2h-2}\overline{Q}_{k,k-h+1}(x^2q^4;q)].
\end{align*}
By \eqref{eqJ}, we have
\begin{align*}
U^{2a}_{2k,2a}(x;q)&=(xq)^{2a+1}(-xq^3;q^2)_\infty\sum_{h=1}^{k-a}[\overline{Q}_{k,h}(x^2q^2;q^2)-\overline{Q}_{k,h-1}(x^2q^2;q^2)]
\\&+(xq)^{2a-1}(-xq^3;q^2)_\infty\sum_{h=1}^{k-a+1}[\overline{Q}_{k,h}(x^2q^2;q^2)-\overline{Q}_{k,h-1}(x^2q^2;q^2)]
\\&=(xq)^{2a+1}(-xq^3;q^2)_\infty\overline{Q}_{k,k-a}(x^2q^2;q^2)+(xq)^{2a-1}(-xq^3;q)_\infty\overline{Q}_{k,k-a+1}(x^2q^2;q^2)
\\&=xq(-xq^3;q^2)_\infty[\overline{Q}_{k,a}(x^2;q^2)-\overline{Q}_{k,a-1}(x^2;q^2)].
\end{align*}
This completes the proof.\qed

\begin{thm}For $k\geq a\geq 1$, we have \begin{align}\nonumber\sum_{n\geq0}U_{2k,2a-1}(n)q^n=&\frac{(-q^2;q)_\infty(q^{2a},q^{4k-2a},q^{4k};q^{4k})_\infty}{(q^2;q^2)_\infty}
\\&+\frac{xq(-q^2;q)_\infty(q^{2a-2},q^{4k-2a+2},q^{4k};q^{4k})_\infty}{(q^2;q^2)_\infty}
\end{align}
\end{thm}
\pf By Theorem \ref{UQ}, \eqref{UU12} and \eqref{U2a}, we compute  $U_{2k,2a-1}(x;q)$ as
\begin{align}\nonumber&U_{2k,2a-1}(x;q)=U_{2k,2a}(x;q)-U^{2a}_{2k,2a}(x;q)
\\\nonumber&=(-xq;q^2)_\infty\overline{Q}_{k,a}(x^2;q^2)-xq(-xq^3;q^2)_\infty\overline{Q}_{k,a}(x^2;q^2)+xq(-xq^3;q^2)_\infty\overline{Q}_{k,a-1}(x^2;q^2)
\\&=(-xq^3;q^2)_\infty\overline{Q}_{k,a}(x^2;q^2)+xq(-xq^3;q^2)_\infty\overline{Q}_{k,a-1}(x^2;q^2).
\end{align}
Let $x=1$, we get the generating function of $U_{2k,2a-1}(n)$.
\begin{align}\nonumber&\sum_{n\geq0}U_{2k,2a-1}(n)q^n=U_{2k,2a-1}(1;q)\\&=\frac{(-q^2;q)_\infty(q^{2a},q^{4k-2a},q^{4k};q^{4k})_\infty}{(q^2;q^2)_\infty}
+\frac{xq(-q^2;q)_\infty(q^{2a-2},q^{4k-2a+2},q^{4k};q^{4k})_\infty}{(q^2;q^2)_\infty}.
\end{align}
\qed

Now we begin to consider the generating function of $U_{2k+1,2a}(n)$. The way is similar as deriving the generating function of $U_{2k,2a-1}(n)$.
\begin{thm}For $k\geq a\geq 1$, we have
\label{12}\begin{equation}\label{UU21}U_{2k+1,2a}(x;q)=U_{2k+1,2a+1}(x;q)-U^{2a+1}_{2k+1,2a+1}(x;q)
\end{equation}
and \begin{equation}\label{U2a}U^{2a}_{2k,2a}(x;q)=xq(-xq^3;q^2)_\infty[\overline{Q}_{k,a}(x^2;q^2)-\overline{Q}_{k,a-1}(x^2;q^2)].
\end{equation}
\end{thm}
The proof is similar as that of Theorem \ref{21}, so we omit it here.

\begin{thm}\begin{align}\nonumber\sum_{n\geq0}U_{2k+1,2a}(n)q^n
=&\frac{(-q^2;q)_\infty(q^{2a+1},q^{4k-2a+1},q^{4k+2};q^{4k+2})_\infty}{(q^2;q^2)_\infty}
\\&+\frac{xq(-q^2;q)_\infty(q^{2a-1},q^{4k-2a+3},q^{4k+2};q^{4k+2})_\infty}{(q^2;q^2)_\infty}.
\end{align}
\end{thm}
\pf By using Theorem \ref{12}, we get
\begin{align}\nonumber &U_{2k+1,2a}(x;q)=U_{2k+1,2a+1}(x;q)-U^{2a+1}_{2k+1,2a+1}(x;q)
\\\nonumber&=(-xq;q^2)_\infty\overline{Q}_{k+\frac{1}{2},a+\frac{1}{2}}(x^2;q^2)
-xq(-xq^3;q^2)_\infty[\overline{Q}_{k+\frac{1}{2},a+\frac{1}{2}}(x^2;q^2)-\overline{Q}_{k+\frac{1}{2},a-\frac{1}{2}}(x^2;q^2)]
\\&=(-xq^3;q^2)_\infty\overline{Q}_{k+\frac{1}{2},a+\frac{1}{2}}(x^2;q^2)+xq(-xq^3;q^2)_\infty\overline{Q}_{k+\frac{1}{2},a-\frac{1}{2}}(x^2;q^2).
\end{align}
Let $x=1$, we reach that
\begin{align}\nonumber&\sum_{n\geq0}U_{2k+1,2a}(n)q^n=U_{2k+1,2a}(1;q)
\\&=\frac{(-q^2;q)_\infty(q^{2a+1},q^{4k-2a+1},q^{4k+2};q^{4k+2})_\infty}{(q^2;q^2)_\infty}
+\frac{xq(-q^2;q)_\infty(q^{2a-1},q^{4k-2a+3},q^{4k+2};q^{4k+2})_\infty}{(q^2;q^2)_\infty}.
\end{align}
\qed

\section{The generating function of $\overline{U}_{2k,2a}(n)$}

Now we begin to derive the generating function of $\overline{U}_{k,a}(n)$. We can derive the following relation by noticing the property that $f_1+f_{\overline{1}}$ is even.
\begin{thm}
For $k\geq a\geq2$ with $a$ even, we have
\begin{equation}\overline{U}_{k,a}(n)=\overline{U}_{k,a-1}(n).
\end{equation}
\end{thm}
In this section, we just derive the generating function of $\overline{U}_{2k,2a}(n)$. By the combinatorial interpretation of $\overline{U}_{2k,2a}(n)$, we can get the relation between $\overline{U}_{2k,2a}(n)$ and $U_{2k,2a}(n)$ as follows.
\begin{thm}For $k\geq a\geq1$, we have\begin{align}\nonumber &\overline{U}_{2k,2a}(x;q)-\overline{U}_{2k,2a-2}(x;q)
\\\label{difU22}&=(xq)^{2a}U_{2k,2k-2a}(xq;q)+(xq)^{2a-2}U_{2k,2k-2a+2}(xq;q).\end{align}
\end{thm}
\pf
By the definition of the overpartitions enumerated by $\overline{U}_{2k,2a}(n)$ and $\overline{U}_{2k,2a-2}(n)$, we know that $\overline{U}_{2k,2a}(n)-\overline{U}_{2k,2a-2}(n)$ is the number of overpartitions that enumerated by $\overline{U}_{2k,2a}(n)$ with the restriction that: if there is an $\overline{1}$, then the number of $1$ is  $2a-1$, otherwise the number of $1$ is  $2a-2$.

Then we consider the overpartition $\lambda$ enumerated by $\overline{U}_{2k,2a}(n)-\overline{U}_{2k,2a-2}(n)$ in the following cases.

(i) $f_{\overline{1}}(\lambda)=1$  and $f_{1}(\lambda)=2a-1$.

 If $f_{\overline{2}}(\lambda)=1$, then $f_2(\lambda)\leq 2k-2a$.
By deleting all the $1$'s, and then subtracting $1$ from each other parts, we get the overpartitions enumerated by $U_{2k,2k-2a}(m-2a,n-m)$ which have the $\overline{1}$ as a part.

If $f_{\overline{2}}(\lambda)=0$, then $f_2(\lambda)\leq 2k-2a-1$.
By deleting all the $1$'s, and then subtracting $1$ from each other parts, we get the overpartitions enumerated by $U_{2k,2k-2a}(m-2a,n-m)$ which do not have the $\overline{1}$ as a part.

(ii) $f_{\overline{1}}(\lambda)=0$  and $f_{1}(\lambda)=2a-2$.

If $f_{\overline{2}}(\lambda)=1$, then $f_2(\lambda)\leq 2k-2a+2$.
By deleting all the $1$'s, and then subtracting $1$ from each other parts, we get the overpartitions enumerated by $U_{2k,2k-2a+2}(m-2a+2,n-m)$ which have the $\overline{1}$ as a part.

If $f_{\overline{2}}(\lambda)=0$, then $f_2(\lambda)\leq 2k-2a+1$.
By deleting all the $1$'s, and then subtracting $1$ from each other parts, we get the overpartitions enumerated by $U_{2k,2k-2a-2}(m-2a+2,n-m)$ which do not have the $\overline{1}$ as a part.
\qed
In the previous section, we have derived the formula of $U_{2k,2a}(x;q)$. Here, we derive the formula of $\overline{U}_{2k,2a}(x;q)$.

\begin{thm}For $k\geq a\geq1$, we have \begin{equation}\sum_{n\geq0}\overline{U}_{2k,2a}(x;q)=(-xq^2;q^2)_\infty\overline{Q}_{k,a}(x^2;q^2).\end{equation}
\end{thm}

\pf We employ Theorem \ref{difU22} and \ref{UQ}, get the following result.
\begin{align*}&\overline{U}_{2k,2a}(x;q)-\overline{U}_{2k,2a-2}(x;q)
\\&=(xq)^{2a}U_{2k,2k-2a}(xq;q)+(xq)^{2a-2}U_{2k,2k-2a+2}(xq;q)
\\&=(xq)^{2a}(-xq^2;q^2)_\infty \overline{Q}_{k,k-a}(x^2q^2;q^2)+(xq)^{2a-2}(-xq^2;q^2)_\infty \overline{Q}_{k,k-a+1}(x^2q^2;q^2).\end{align*}
By using \eqref{eqJ}, we have
\begin{equation}\label{overlineUQ}\overline{U}_{2k,2a}(x;q)-\overline{U}_{2k,2a-2}(x;q)
 =(-xq^2;q^2)_\infty(\overline{Q}_{k,a}(x^2;q^2)-\overline{Q}_{k,a-1}(x^2;q^2)).
\end{equation}
 We successively  apply  \eqref{overlineUQ} to get the following relation

\begin{equation}\overline{U}_{2k,2a}(x;q) =(-xq^2;q^2)_\infty\overline{Q}_{k,a}(x^2;q^2).
\end{equation}
\qed
  Letting $x=1$, we get the generating function of  $\overline{U}_{2k,2a}(n)$.
 \begin{thm}For $k\geq a\geq 1$, we have
 \begin{equation}\sum_{n\geq 0}\overline{U}_{2k,2a}(n)q^n=\frac{(-q^2;q^2)^2_\infty(q^{2a},q^{4k-2a},q^{4k};q^{4k})_\infty}{(q^{2};q^2)_\infty}.
 \end{equation}
 \end{thm}
\section{The generating function of $U_{2k+1,2a+1}(n)$ and   $\overline{U}_{2k+1,2a+1}(n)$}

In this section we consider the generating function $U_{2k+1,2a+1}(x;q)$ and   $\overline{U}_{2k+1,2a+1}(x;q)$. By using the combinatorial tool we give
 the $q$-differential relations between   $U_{2k+1,2a+1}(x;q)$ and   $\overline{U}_{2k+1,2a+1}(x;q)$. Then by the initial values, we can get the   formulas of
 $U_{2k+1,2a+1}(x;q)$ and   $\overline{U}_{2k+1,2a+1}(x;q)$.
We first check the following relations.
\begin{thm}\label{Uoo}For $k\geq a\geq1$, we have \begin{align}\nonumber&U_{2k+1,2a+1}(x;q)-U_{2k+1,2a-1}(x;q)
\\\label{difU11}&=(1+xq)(xq)^{2a+1}\overline{U}_{2k+1,2k-2a}(xq;q)+(1+xq)(xq)^{2a-1}\overline{U}_{2k+1,2k-2a+2}(xq;q),\end{align}
 and
  \begin{equation}\label{difU1}U_{2k+1,1}(x;q)=x^2q^2\overline{U}_{2k+1,2k}(xq;q)+\overline{U}_{2k+1,2k+2}(xq;q).
  \end{equation}

\end{thm}
\pf We prove the relation \eqref{difU11} by consider the overpartitions enumerated by $U_{2k+1,2a+1}(m,n)-U_{2k+1,2a-1}(m,n)$.
\begin{itemize}
\item[1.]there is an $\overline{1}$ and there are $2a+1$  or $2a$ $1$'s. Since $f_2+f_{\overline{2}}$ is even,  if there is an $\overline{2}$, then there are at most $2k-2a-1$ non-overlined $2$'s or if there is no $\overline{2}$,  then there are at most $2k-2a-2$ nonoverlined $2$'s.
    By deleting all $1$'s and $\overline{1}$ and subtracting $1$ from each of other parts, we get all overpartitions enumerated by $\overline{U}_{2k+1,2k-2a}(m-2a-1,n-m)$ and $\overline{U}_{2k+1,2k-2a}(m-2a-2,n-m)$. Then the generating function is $(1+xq)(xq)^{2a+1}\overline{U}_{2k+1,2k-2a}(xq;q)$.
\item[2.]There is no $\overline{1}$ and there are $2a$ or $2a-1$ $1$'s.  since $f_2+f_{\overline{2}}$ is even, if there is an $\overline{2}$, then there are at most $2k-2a+1$ nonoverlined $2$'s or if there is no $\overline{2}$ there are at most $2k-2a$ nonoverlined $2$'s.
    By deleting all $1$'s and subtracting $1$ from each other parts, we get all overpartitions enumerated by $\overline{U}_{2k+1,2k-2a+2}(m-2a,n-m)$ and $\overline{U}_{2k+1,2k-2a+1}(m-2a+2,n-m)$. Then the generating function is $(1+xq)(xq)^{2a-1}\overline{U}_{2k+1,2k-2a+2}(xq;q)$.
\end{itemize}

 To prove \eqref{difU1}, we consider the overpartitions enumerated by $U_{2k+1,1}(m,n)$. If there is an $\overline{1}$ then there is a nonoverlined $1$, otherwise there are no $1$'s.  Deleting all $1$ and $\overline{1}$ and subtracting $1$ from each other part, we get the overpartitions enumerated by $\overline{U}_{2k+1,2k-1}(m-2,n-m)$ and the overpartitions enumerated by $\overline{U}_{2k+1,2a+1}(m,n-m)$. The generating function is $x^2q^2\overline{U}_{2k+1,2k-1}(xq;q)+\overline{U}_{2k+1,2k+1}(xq;q)$. This completes the proof.\qed

\begin{thm}For $k\geq a\geq1$, we have
\begin{align}
\nonumber&\overline{U}_{2k+1,2a}(n)-\overline{U}_{2k+1,2a-2}(n)
\\&=(xq)^{2a}U_{2k+1,2k-2a+1}(xq;q)+(xq)^{2a-2}U_{2k+1,2k-2a+3}(xq;q),
\end{align}
and \begin{equation}\overline{U}_{2k+1,2}(x;q)=(xq)^2U_{2k+1,2k-1}(xq;q)+U_{2k+1,2k+1}(xq;q).\end{equation}
\end{thm}
The proof is similar as Theorem \ref{Uoo}, so we omit it here.

Now we derive the following result, by which we can get the generating function of $U_{2k+1,2a+1}(n)$ and $\overline{U}_{2k+1,2a}(n)$.
\begin{thm}For $k\geq a\geq0$, we have 
\begin{equation}U_{2k+1,2a+1}(x;q)=(-xq;q^2)_\infty\overline{Q}_{k+\frac{1}{2},a+\frac{1}{2}}(x^2;q^2).
\end{equation}
\begin{equation}\overline{U}_{2k+1,2a}(x;q)=\overline{U}_{2k+1,2a-1}(x;q)=(-xq^2;q^2)_\infty\overline{Q}_{k+\frac{1}{2},a}(x^2;q^2).
\end{equation}
\end{thm}

\pf
Let \[T_{2k+1,2a+1}(x;q):=(-xq;q^2)_\infty\overline{Q}_{k+\frac{1}{2},a+\frac{1}{2}}(x^2;q^2)\]
and  \begin{equation}\label{TT12}\overline{T}_{2k+1,2a}(x;q):=(-xq^2;q^2)_\infty\overline{Q}_{k+\frac{1}{2},a}(x^2;q^2).\end{equation}
We shall prove that
\[T_{2k+1,2a+1}(x;q)=U_{2k+1,2a+1}(x;q)\]
and \[\overline{T}_{2k+1,2a}(x;q)=\overline{U}_{2k+1,2a}(x;q)=\overline{U}_{2k+1,2a-1}(x;q).\]
We first check the initial values

 One can easily derive the following result by Lemma \ref{lemQini}
\begin{equation}
1 =T_{2k+1,2a+1}(0; q) = T_{2k+1,2a+1}(x; 0) = \overline{T}_{2k+1,2a}(0; q) = \overline{T}_{2k+1,2a}(x; 0).
\end{equation}
By \eqref{Qminus} and Lemma \ref{lemQ}, one can derive the following relation
 \begin{align*}&T_{2k+1,1}(x;q)
 =(-xq;q^2)_\infty\overline{Q}_{k+\frac{1}{2},\frac{1}{2}}(x^2;q^2)
 \\&=(-xq^3;q^2)_\infty(1+xq)\overline{Q}_{k+\frac{1}{2},\frac{1}{2}}(x^2;q^2)
 \\&=(-xq^3;q^2)_\infty xq[\overline{Q}_{k+\frac{1}{2},\frac{1}{2}}(x^2;q^2)-\overline{Q}_{k+\frac{1}{2},-\frac{1}{2}}(x^2;q^2)]
 \\&= (-xq^3;q^2)_\infty xq[(xq)\overline{Q}_{k+\frac{1}{2},k}(x^2q^2;q^2)-(xq)^{-1}\overline{Q}_{k+\frac{1}{2},k+1}(x^2q^2;q^2)]
 \\&=x^2q^2\overline{T}_{2k+1,2k}(xq;q)+\overline{T}_{2k+1,2k+2}(xq;q).
 \end{align*}
One also can check that
\begin{align}\nonumber&T_{2k+1,2a+1}(x;q)-T_{2k+1,2a-1}(x;q)
\\\nonumber&=(-xq;q^2)_\infty[\overline{Q}_{k+\frac{1}{2},a+\frac{1}{2}}(x^2;q^2)-\overline{Q}_{k+\frac{1}{2},a-\frac{1}{2}}(x^2;q^2)]
\\\nonumber&=(-xq;q^2)_\infty[(xq)^{2a+1}\overline{Q}_{k+\frac{1}{2},k-a}(x^2q^2;q^2)+(xq)^{2a-1}\overline{Q}_{k+\frac{1}{2},k-a+1}(x^2q^2;q^2)].
\end{align}
By \eqref{TT12},     \[\overline{T}_{2k+1,2a}(xq;q)=(-xq^3;q^2)_\infty\overline{Q}_{k+\frac{1}{2},a}(x^2q^2;q^2),\]
then we have
\begin{align}\nonumber&T_{2k+1,2a+1}(x;q)-T_{2k+1,2a-1}(x;q)\\&=(1+xq)[(xq)^{2a+1}\overline{T}_{2k+1,2k-2a}(xq;q)+(xq)^{2a-1}\overline{T}_{2k+1,2k-2a+2}(xq;q)].
\end{align}
 Similarly, we can get the following relation
\begin{align}\nonumber&\overline{T}_{2k+1,2a}(x;q)-\overline{T}_{2k+1,2a-2}(x;q)
\\\nonumber&=(-xq^2;q^2)_\infty[\overline{Q}_{k+\frac{1}{2},a}(x^2;q^2)-\overline{Q}_{k+\frac{1}{2},a-1}(x^2;q^2)]
\\\nonumber&=(-xq^2;q^2)_\infty[(x^2q^2)^{a}\overline{Q}_{k+\frac{1}{2},k-a+\frac{1}{2}}(x^2q^2;q^2)
+(x^2q^2)^{a-1}\overline{Q}_{k+\frac{1}{2},k-a+\frac{3}{2}}(x^2q^2;q^2)]
\\&=(xq)^{2a}T_{2k+1,2k-2a+1}(xq;q)+(xq)^{2a-2}T_{2k+1,2k-2a+3}(xq;q). \end{align}
Now we consider the initial value of $\overline{T}_{2k+1,2}(x;q)$.
\begin{align}
&\nonumber\overline{T}_{2k+1,2}(x;q)=(-xq^2;q^2)\overline{Q}_{k+\frac{1}{2},1}(x^2;q^2)
\\\nonumber=&(-xq^2;q^2)[\overline{Q}_{k+\frac{1}{2},1}(x^2;q^2)-\overline{Q}_{k+\frac{1}{2},0}(x^2;q^2)]
\\\nonumber&=(-xq^2;q^2)(x^2q^2)\overline{Q}_{k+\frac{1}{2},k-\frac{1}{2}}(x^2q^2)+\overline{Q}_{k+\frac{1}{2},k-\frac{1}{2}}(x^2;q^2)
\\&=(xq)^2T_{2k+1,2k-1}(xq;q)+T_{2k+1,2k+1}(xq;q)
\end{align}
By the initial values and double inductions we can derive that
\[T_{2k+1,2a+1}(x;q)=U_{2k+1,2a+1}(x;q)\]
and \[\overline{T}_{2k+1,2a}(x;q)=\overline{U}_{2k+1,2a}(x;q)=\overline{U}_{2k+1,2a-1}(x;q).\]
Then we completes the proof.

Letting $x=1$, we get the following theorem.
\begin{thm}
\begin{equation}
\sum_{n=0}^\infty U_{2k+1,2a+1}(n)q^n=\frac{(-q;q)_\infty(q^{2a+1},q^{4k+1-2a},q^{4k+2};q^{4k+2})_\infty}{(q^2;q^2)_\infty}.
\end{equation}
\begin{equation}\sum_{n=0}^\infty \overline{U}_{2k+1,2a}(n)q^n=\sum_{n=0}^\infty \overline{U}_{2k+1,2a-1}(n)q^n=\frac{(-q^2;q^2)^2_\infty(q^{2a},q^{4k-2a+2},q^{4k+2};q^{4k+2})_\infty}{(q^2;q^2)_\infty}.
\end{equation}
\end{thm}

\section{The combinatorial relations with $\overline{B}_{k,a}(n)$}

By the generating function, we can get some  relations among $U_{k,a}(n)$, $\overline{U}_{k,a}(n)$ and $\overline{B}_{k,a}(n)$. In this section, we shall give the combinatorial interpretation of these relations.
\begin{thm}For $k\geq a\geq1$, we have 
\begin{equation}\label{UB}\sum_{n\geq0}U_{2k,2a}(n)q^n=(-q;q^2)_\infty\sum_{n\geq 0}\overline{B}_{k,a}(n)q^{2n},\end{equation}
and \begin{equation}\label{UB2}\sum_{n\geq0}\overline{U}_{2k,2a}(n)q^n=(-q^2;q^2)_\infty\sum_{n\geq 0}\overline{B}_{k,a}(n)q^{2n}.
\end{equation}
\end{thm}
\pf Firstly, we prove \eqref{UB}. Let $\lambda=\lambda_1+\lambda_2+\cdots+\lambda_m$ be a partition enumerated by $U_{2k,2a}(n)$, i.e., $f_l(\lambda)+f_{\overline{l}}(\lambda)+f_{l+1}(\lambda)\leq 2-1+f_{\overline{l+1}}(\lambda)$ and $f_{2l}(\lambda)+f_{\overline{2l}}(\lambda)$ are even for each $l$.  

If there is a number $2l-1$ such that $f_{2l-1}(\lambda)+f_{\overline{2l-1}}(\lambda)$ is odd, then we remove one nonoverlined $2l-1$.
 We let
$\gamma$ be the partition consisting of those removed odd parts and $\beta$ be  the overpartition of the remaining parts
of $\lambda$. Then it is obvious that $\gamma$ has distinct odd parts, and $\beta$ is an overpartition such that $f_{l}(\beta)+f_{\overline{l}}(\beta)$ are even, for each $l$ with $f_1(\beta)\leq 2a-1+f_{\overline{1}}(\beta)$.

Now, we add two repeated parts of $\beta$ and then divide them by $2$,  to obtain an overpartition $\beta'$. If there is an overlined $l$ in $\beta$, then there is an $\overline{l}$ in $\beta'$.  We will check that $\beta'$ is an overpartition enumerated by $\overline{B}_{k,a}((n-|\gamma|)/2)$.

If $f_{\overline{l}}(\lambda)+f_{l}(\lambda)+f_{\overline{l+1}}(\lambda)+f_{l+1}(\lambda)\leq 2k-2$, then after the operation, $f_{\overline{l}}(\beta')+f_{l}(\beta')+f_{\overline{l+1}}(\beta')+f_{l+1}(\beta')\leq k-1$.

If $f_{\overline{l}}(\lambda)+f_{l}(\lambda)+f_{\overline{l+1}}(\lambda)+f)_{l+1}(\lambda)=2k-1$, then the odd number in $\{l,l+1\}$ appears odd number of times (including overlined and nonoverlined) so that  it should be removed. Thus, after the operation, $f_{\overline{l}}(\beta')+f_{l}(\beta')+f_{\overline{l+1}}(\beta')+f_{l+1}(\beta')= k-1$.

If $f_{\overline{l}}(\lambda)+f_{l}(\lambda)+f_{\overline{l+1}}(\lambda)+f_{l+1}(\lambda)=2k$ then $f_{\overline{l+1}}(\lambda)=1$. Hence, after the operation, $f_{\overline{l}}(\beta')+f_{l}(\beta')+f_{\overline{l+1}}(\beta')+f_{l+1}(\beta')= k$ and $f_{\overline{l+1}}(\beta')=1$, so $f_{\overline{l}}(\beta')+f_{l}(\beta')+f_{l+1}(\beta')= k-1$.

Then we can see that $\beta'$ is an overpartition that $f_{\overline{l}}(\beta')+f_{l}(\beta')+f_{l+1}(\beta')\leq k-1$.

Now we check that $f_1(\beta')\leq a-1$. If $f_{\overline{1}}(\lambda)=1$ then together with that $f_1(\lambda)+f_{\overline{1}}(\lambda)\leq 2a$, we have $f_1(\beta')+f_{\overline{1}}(\beta')\leq a$ and  $f_{\overline{1}}(\beta')=1$.  If $f_{\overline{1}}(\lambda)=0$, then $f_1\leq 2a-1$. If $f_1(\lambda)\neq2a-1$, then
$f_1(\beta')\leq a-1$, otherwise  $f_1(\lambda)$ is odd, so we should remove  one $1$. Then after the operation $f_1(\beta')=a-1$.

Now, we have proved that $\beta'$ is an overpartition enumerated by $\overline{B}_{k,a}((n-|\gamma|)/2)$, i.e. \eqref{UB}.

The combinatorial interpretation of \eqref{UB2} is similar, so we omit it here. \qed

\vspace{0.5cm}
 \noindent{\bf Acknowledgments.} This work was supported by the National Science Foundation of China (Nos.1140149, 11501089, 11501408).

\end{document}